\magnification\magstephalf
\documentstyle{amsppt}

\hsize 5.72 truein
\vsize 7.9 truein
\hoffset .39 truein
\voffset .26 truein
\mathsurround 1.67pt
\parindent 20pt
\normalbaselineskip 13.8truept
\normalbaselines
\binoppenalty 10000
\relpenalty 10000
\csname nologo\endcsname 


\font\bc=cmb10
\font\tenbsy=cmbsy10

\catcode`\@=11

\def\myitem#1.{\item"(#1)."\advance\leftskip10pt\ignorespaces}

\def\qedsymbol{{\mathsurround\z@$\square$}}
\redefine\qed{\relaxnext@\ifmmode\let\next\@qed\else
  {\unskip\nobreak\hfil\penalty50\hskip2em\null\nobreak\hfil
    \qedsymbol\parfillskip\z@\finalhyphendemerits0\par}\fi\next}
\def\@qed#1$${\belowdisplayskip\z@\belowdisplayshortskip\z@
  \postdisplaypenalty\@M\relax#1
  $$\par{\lineskip\z@\baselineskip\z@\vbox to\z@{\vss\noindent\qed}}}
\outer\redefine\beginsection#1#2\par{\par\penalty-250\bigskip\vskip\parskip
  \leftline{\tenbsy x\bf#1. #2}\nobreak\smallskip\noindent}
\outer\redefine\genbeginsect#1\par{\par\penalty-250\bigskip\vskip\parskip
  \leftline{\bf#1}\nobreak\smallskip\noindent}

\def\next{\let\@sptoken= }\def\next@{ }\expandafter\next\next@
\def\@futureletnext#1{\let\nextii@#1\futurelet\next\@flti}
\def\@flti{\ifx\next\@sptoken\let\next@\@fltii\else\let\next@\nextii@\fi\next@}
\expandafter\def\expandafter\@fltii\next@{\futurelet\next\@flti}

\let\zeroindent\z@
\let\savedef@\endproclaim\let\endproclaim\relax 
\define\chkproclaim@{\add@missing\endroster\add@missing\enddefinition
  \add@missing\endproclaim
  \envir@stack\endproclaim
  \edef\endit@{\leftskip\the\leftskip\rightskip\the\rightskip}}
\let\endproclaim\savedef@
\def\thing@{.\enspace\egroup\ignorespaces}
\def\thingi@(#1){ \rm(#1)\thing@}
\def\thingii@\cite#1{ \rm\@pcite{#1}\thing@}
\def\thingiii@{\ifx\next(\let\next\thingi@
  \else\ifx\next\cite\let\next\thingii@\else\let\next\thing@\fi\fi\next}
\def\thing#1#2#3{\chkproclaim@
  \ifvmode \medbreak \else \par\nobreak\smallskip \fi
  \noindent\advance\leftskip#1
  \hskip-#1#3\bgroup\bc#2\unskip\@futureletnext\thingiii@}
\let\savedef@\endproclaim\let\endproclaim\relax 
\def\endit{\endproclaim\endit@\let\endit@\undefined}
\let\endproclaim\savedef@
\def\defn#1{\thing\parindent{Definition #1}\rm}
\def\lemma#1{\thing\parindent{Lemma #1}\sl}
\def\prop#1{\thing\parindent{Proposition #1}\sl}

\def\conj#1{\thing\parindent{Conjecture #1}\sl}

\def\remk#1{\thing\zeroindent{Remark #1}\rm}
\def\example#1{\thing\zeroindent{Example #1}\rm}
\def\narrowthing#1{\chkproclaim@\medbreak\narrower\noindent
  \it\def\next{#1}\def\next@{}\ifx\next\next@\ignorespaces
  \else\bgroup\bc#1\unskip\let\next\narrowthing@\fi\next}
\def\narrowthing@{\@futureletnext\thingiii@}

\def\@cite#1,#2\end@{{\rm([\bf#1\rm],#2)}}
\def\cite#1{\in@,{#1}\ifin@\def\next{\@cite#1\end@}\else
  \relaxnext@{\rm[\bf#1\rm]}\fi\next}
\def\@pcite#1{\in@,{#1}\ifin@\def\next{\@cite#1\end@}\else
  \relaxnext@{\rm([\bf#1\rm])}\fi\next}

\advance\minaw@ 1.2\ex@
\atdef@[#1]{\ampersand@\let\@hook0\let\@twohead0\brack@i#1,\z@,}
\def\brack@{\z@}
\let\@@hook\brack@
\let\@@twohead\brack@
\def\brack@i#1,{\def\next{#1}\ifx\next\brack@
  \let\next\brack@ii
  \else \expandafter\ifx\csname @@#1\endcsname\brack@
    \expandafter\let\csname @#1\endcsname1\let\next\brack@i
    \else \Err@{Unrecognized option in @[}%
  \fi\fi\next}
\def\brack@ii{\futurelet\next\brack@iii}
\def\brack@iii{\ifx\next>\let\next\brack@gtr
  \else\ifx\next<\let\next\brack@less
    \else\relaxnext@\Err@{Only < or > may be used here}
  \fi\fi\next}
\def\brack@gtr>#1>#2>{\setboxz@h{$\m@th\ssize\;{#1}\;\;$}%
 \setbox@ne\hbox{$\m@th\ssize\;{#2}\;\;$}\setbox\tw@\hbox{$\m@th#2$}%
 \ifCD@\global\bigaw@\minCDaw@\else\global\bigaw@\minaw@\fi
 \ifdim\wdz@>\bigaw@\global\bigaw@\wdz@\fi
 \ifdim\wd@ne>\bigaw@\global\bigaw@\wd@ne\fi
 \ifCD@\enskip\fi
 \mathrel{\mathop{\hbox to\bigaw@{$\ifx\@hook1\lhook\mathrel{\mkern-9mu}\fi
  \setboxz@h{$\displaystyle-\m@th$}\ht\z@\z@
  \displaystyle\m@th\copy\z@\mkern-6mu\cleaders
  \hbox{$\displaystyle\mkern-2mu\box\z@\mkern-2mu$}\hfill
  \mkern-6mu\mathord\ifx\@twohead1\twoheadrightarrow\else\rightarrow\fi$}}%
 \ifdim\wd\tw@>\z@\limits^{#1}_{#2}\else\limits^{#1}\fi}%
 \ifCD@\enskip\fi\ampersand@}
\def\brack@less<#1<#2<{\setboxz@h{$\m@th\ssize\;\;{#1}\;$}%
 \setbox@ne\hbox{$\m@th\ssize\;\;{#2}\;$}\setbox\tw@\hbox{$\m@th#2$}%
 \ifCD@\global\bigaw@\minCDaw@\else\global\bigaw@\minaw@\fi
 \ifdim\wdz@>\bigaw@\global\bigaw@\wdz@\fi
 \ifdim\wd@ne>\bigaw@\global\bigaw@\wd@ne\fi
 \ifCD@\enskip\fi
 \mathrel{\mathop{\hbox to\bigaw@{$%
  \setboxz@h{$\displaystyle-\m@th$}\ht\z@\z@
  \displaystyle\m@th\mathord\ifx\@twohead1\twoheadleftarrow\else\leftarrow\fi
  \mkern-6mu\cleaders
  \hbox{$\displaystyle\mkern-2mu\copy\z@\mkern-2mu$}\hfill
  \mkern-6mu\box\z@\ifx\@hook1\mkern-9mu\rhook\fi$}}%
 \ifdim\wd\tw@>\z@\limits^{#1}_{#2}\else\limits^{#1}\fi}%
 \ifCD@\enskip\fi\ampersand@}


\define\today{\number\day\ \ifcase\month\or
  January\or February\or March\or April\or May\or June\or
  July\or August\or September\or October\or November\or December\fi
  \ \number\year}
\def\pr@m@s{\ifx'\next\let\nxt\pr@@@s \else\ifx^\next\let\nxt\pr@@@t
  \else\let\nxt\egroup\fi\fi \nxt}

\define\widebar#1{\mathchoice
  {\setbox0\hbox{\mathsurround\z@$\displaystyle{#1}$}\dimen@.1\wd\z@
    \ifdim\wd\z@<.4em\relax \dimen@ -.16em\advance\dimen@.5\wd\z@ \fi
    \ifdim\wd\z@>2.5em\relax \dimen@.25em\relax \fi
    \kern\dimen@ \overline{\kern-\dimen@ \box0\kern-\dimen@}\kern\dimen@}%
  {\setbox0\hbox{\mathsurround\z@$\textstyle{#1}$}\dimen@.1\wd\z@
    \ifdim\wd\z@<.4em\relax \dimen@ -.16em\advance\dimen@.5\wd\z@ \fi
    \ifdim\wd\z@>2.5em\relax \dimen@.25em\relax \fi
    \kern\dimen@ \overline{\kern-\dimen@ \box0\kern-\dimen@}\kern\dimen@}%
  {\setbox0\hbox{\mathsurround\z@$\scriptstyle{#1}$}\dimen@.1\wd\z@
    \ifdim\wd\z@<.28em\relax \dimen@ -.112em\advance\dimen@.5\wd\z@ \fi
    \ifdim\wd\z@>1.75em\relax \dimen@.175em\relax \fi
    \kern\dimen@ \overline{\kern-\dimen@ \box0\kern-\dimen@}\kern\dimen@}%
  {\setbox0\hbox{\mathsurround\z@$\scriptscriptstyle{#1}$}\dimen@.1\wd\z@
    \ifdim\wd\z@<.2em\relax \dimen@ -.08em\advance\dimen@.5\wd\z@ \fi
    \ifdim\wd\z@>1.25em\relax \dimen@.125em\relax \fi
    \kern\dimen@ \overline{\kern-\dimen@ \box0\kern-\dimen@}\kern\dimen@}%
  }

\catcode`\@\active

\let\PVstyle=d 

\loadeufm

\font\tenscr=rsfs10 
\font\sevenscr=rsfs7 
\font\fivescr=rsfs5 
\skewchar\tenscr='177 \skewchar\sevenscr='177 \skewchar\fivescr='177
\newfam\scrfam \textfont\scrfam=\tenscr \scriptfont\scrfam=\sevenscr
\scriptscriptfont\scrfam=\fivescr
\define\scr#1{{\fam\scrfam#1}}
\let\Cal\scr

\let\0\relax 
\define\exc{{\text{exc}}}

\define\red{{\text{red}}}
\define\restrictedto#1{\big|_{#1}}

\topmatter
\title Multiplier Ideal Sheaves, Nevanlinna Theory, and Diophantine
  Approximation\endtitle
\rightheadtext{multiplier ideal sheaves}
\author Paul Vojta\endauthor
\affil University of California, Berkeley\endaffil
\address Department of Mathematics, University of California,
  970 Evans Hall\quad\#3840, Berkeley, CA \ 94720-3840\endaddress
\date \today \enddate
\thanks Supported by NSF grants DMS-0200892 and DMS-0500512.\endthanks

\abstract
This note states a conjecture for Nevanlinna theory or diophantine
approximation, with a sheaf of ideals in place of the normal crossings divisor.
This is done by using a correction term involving a multiplier ideal sheaf.
This new conjecture trivially implies earlier conjectures in Nevanlinna theory
or diophantine approximation, and in fact is equivalent to these conjectures.
Although it does not provide anything new, it may be a more convenient
formulation for some applications.
\endabstract
\endtopmatter

\document

This note states a conjecture for Nevanlinna theory or
diophantine approximation, with a sheaf of ideals in place of the normal
crossings divisor.  This is done by using a correction term involving a
multiplier ideal sheaf.  This new conjecture
is equivalent to earlier conjectures in Nevanlinna theory or diophantine
approximation, but may be a more convenient formulation for some applications.
It also shows how multiplier ideal sheaves may have a role in Nevanlinna theory
and diophantine approximation, and therefore may give more information on
the structure of the situation.

Section \01 briefly describes multiplier ideal sheaves, and gives a
variant definition specific to this situation.  Section \02 describes
proximity functions for sheaves of ideals, using work of Silverman and
Yamanoi.  Sections \03 and \04 form the heart of the paper, giving the
conjectures and showing their equivalence to previous conjectures.

\narrowthing{}
Throughout this paper, $X$ is a smooth complete variety over $\Bbb C$
(in the case of Nevanlinna theory) or over a global field of characteristic
zero (in the case of diophantine approximation).
\endit

\beginsection{\01}{Multiplier Ideal Sheaves}

\defn{\01.1}  Let $\frak a$ be a nonzero sheaf of ideals on $X$, and let
$c\in\Bbb R_{\ge0}$.  Let $\mu\:X'\to X$ be a proper birational morphism
such that $X'$ is a smooth variety and
$$\mu^{*}(\frak a)=\Cal O_{X'}(-F)$$
for a divisor $F$ on $X'$ with normal crossings support.  Let $K_{X'/X}$
denote the ramification divisor of $X'$ over $X$.  Then the {\bc multiplier
ideal sheaf} associated to $\frak a$ and $c$ is the ideal sheaf
$$\Cal I(\frak a^c) = \mu_{*}\Cal O_{X'}\bigl(K_{X'/X}-\lfloor cF\rfloor\bigr)
  \;.$$
\endit

By a theorem of Esnault and Viehweg \cite{L, Thm.~9.2.18}, this definition
is independent of the choice of $\mu$.

For our purposes we need a slightly different definition.

\defn{\01.2}  Let $\frak a$ and $c$ be as above.  We then define
$$\Cal I^{-}(\frak a^c) = \lim_{\epsilon\to0^{+}}\Cal I(\frak a^{c-\epsilon})
  \;.$$
Here we use the discrete topology on the set of ideal sheaves on $X$, and note
that the limit exists because there are only finitely many coefficients
in $\lfloor(c-\epsilon)F\rfloor$.
\endit

We also write $\Cal I(\frak a)=\Cal I(\frak a^1)$ and
$\Cal I^{-}(\frak a)=\Cal I^{-}(\frak a^1)$.

\example{\01.3}  Let $D$ be an (effective, reduced) normal crossings divisor
on $X$ and let $\frak a=\Cal O(-D)$.  Then we can take $X'=X$, in which case
$F=D$ and
$$K_{X'/X}=\lfloor(1-\epsilon)F\rfloor=0\;,$$
so $\Cal I^{-}(\frak a)=\Cal O_X$ (the ideal sheaf corresponding to the
empty closed subscheme).  More generally, if $D$ is effective and has
normal crossings support but is not necessarily reduced, then
$\lfloor(1-\epsilon)F\rfloor=D-D_\red$, and therefore
$\Cal I^{-}(\Cal O(-D))=\Cal O(-(D-D_\red))$.
\endit

\beginsection{\02}{Proximity Functions for Ideal Sheaves}

Silverman \cite{S, 2.2} introduced Weil functions associated to sheaves
of ideals on $X$.  By \cite{S, Thm.~2.1}, there is a unique way to associate
to each ideal sheaf $\frak a\ne(0)$ of $X$ a Weil-like function
$\lambda_{\frak a}$ on $X\setminus Y$, where $Y$ is the closed subscheme
associated to $\frak a$, such that $\lambda_{\frak a}=\lambda_D$ is a
Weil function in the usual sense if $\frak a=\Cal O(-D)$ for some effective
Cartier divisor $D$, and
$\lambda_{\frak a+\frak b}=\min\{\lambda_{\frak a},\lambda_{\frak b}\}$
for all nonzero ideal sheaves $\frak a$ and $\frak b$ of $X$.  Here
uniqueness and equality are up to addition of functions bounded by
$M_k$\snug-constants.

These Weil functions also satisfy the following conditions:
\roster
\item"$\bullet$" They are functorial in the sense that if $f\:X'\to X$
is a morphism of complete varieties with $f(X')\nsubseteq Y$, then
$\lambda_{f^{*}\frak a}=\lambda_{\frak a}\circ f$.
\item"$\bullet$" If $\frak a\subseteq\frak b$ are ideal sheaves on $X$,
then $\lambda_{\frak a}\ge\lambda_{\frak b}$.
\endroster

See also Noguchi \cite{N} and Yamanoi \cite{Y, 2.2}.  They used similar
Weil functions to define proximity functions relative to ideal sheaves.
These are defined as follows.  Let $f\:\Bbb C\to X$ be a holomorphic curve
whose image is not entirely contained in $Y$.  Then we define the proximity
function in the usual way:
$$m_f(\frak a,r)
  = \int_0^{2\pi} \lambda_{\frak a}(f(re^{i\theta}))\frac{d\theta}{2\pi}\;,$$
and similarly in the diophantine case.  Again, this agrees with $m_f(D,r)$
(up to $O(1)$) if $\frak a=\Cal O(-D)$, and satisfies the above two additional
properties (again, up to $O(1)$).

\beginsection{\03}{The Conjectures}

In Nevanlinna theory, we can make the following conjecture:

\conj{\03.1}  Let $X$ be a nonsingular complete complex variety, let $K$ be the
canonical divisor class on $X$, let $\frak a\ne(0)$ be an ideal sheaf on $X$,
let $A$ be a big divisor on $X$, and let $\epsilon>0$.  Then there is a proper
Zariski-closed subset $Z$ of $X$, depending only on $X$, $\frak a$, $A$,
and $\epsilon$, such that if $f\:\Bbb C\to X$ is a holomorphic curve whose
image is not contained in $Z$, then
$$T_{K,f}(r) + m_f(\frak a,r) - m_f(\Cal I^{-}(\frak a),r)
  \le_\exc \epsilon\,T_{A,f}(r) + O(1)\;.$$
Here the subscript ``exc'' means that the inequality holds outside of a
set of $r$ of finite Lebesgue measure.
\endit

One can also make the corresponding conjecture for algebroid functions.

The corresponding conjecture in number theory is:

\conj{\03.2}  Let $k$ be a global field of characteristic zero, let $S$ be
a finite set of places of $k$ containing all archimedean places, let $r$ be a
positive integer, let $X$ be a nonsingular complete variety over $k$,
let $K$ be the canonical divisor class of $X$, let $\frak a\ne(0)$ be an
ideal sheaf on $X$, let $A$ be a big divisor on $X$, and let $\epsilon>0$.
Then there is a proper Zariski-closed subset $Z$ of $X$, depending only on
$k$, $S$, $X$, $\frak a$, $r$, $A$, and $\epsilon$, such that
$$h_K(P) + m(\frak a,P) - m(\Cal I^{-}(\frak a),P)
  \le \epsilon\,h_A(P) + d(P) + O(1)$$
for all $P\in (X\setminus Z)(\bar k)$ with $[k(P):k]\le r$.
\endit

These conjectures obviously generalize earlier conjectures in each case.
Indeed, let $D$ be an effective, reduced normal crossings divisor and
let $\frak a=\Cal O(-D)$.  Then $m_f(\Cal I^{-}(\frak a),r)=O(1)$ and
$m_f(\frak a,r)=m_f(D,r)$, and likewise in the diophantine case.

\prop{\03.3}  Conjectures \03.1 and \03.2 are equivalent to their respective
special cases in which $\frak a=\Cal O(-D)$ with $D$ as above.
\endit

\demo{Proof}  Let $\mu\:X'\to X$, $K_{X'/X}$, and $F$ be as in the definition
of multiplier ideal sheaf, and choose $\eta>0$ such that
$\Cal I^{-}(\frak a)=\Cal I(\frak a^{1-\eta})$.  In the Nevanlinna case,
let $g\:\Bbb C\to X'$ be a lifting of $f$; then
$$\split & T_{K_X,f}(r) + m_f(\frak a,r) - m_f(\Cal I^{-}(\frak a),r) \\
  &\qquad \le T_{K_{X'},g}(r) - m_g(K_{X'/X},r) + m_g(F,r)
    - m_g(-K_{X'/X}+\lfloor(1-\eta)F\rfloor,r) + O(1) \\
  &\qquad = T_{K_{X'},g}(r) + m_g(F_\red,r) + O(1) \\
  &\qquad \le_\exc \epsilon\,T_{A,f}(r) + O(1)\;.\endsplit$$
Here we use the fact that
$$\mu^{*}\Cal I(\frak a^{1-\eta})
  =\mu^{*}\mu_{*}\Cal O_{X'}(K_{X'/X}-\lfloor(1-\eta)F\rfloor)
  \subseteq\Cal O_{X'}(K_{X'/X}-\lfloor(1-\eta)F\rfloor)$$
and therefore
$$\split m_f(\Cal I(\frak a^{1-\eta}),r)
  &\ge m_g(\Cal O_{X'}(K_{X'/X}-\lfloor(1-\eta)F\rfloor),r) + O(1) \\
  &= m_g(-K_{X'/X}+\lfloor(1-\eta)F\rfloor,r) + O(1)\;.\endsplit$$

The diophantine case is similar and is left to the reader.\qed
\enddemo

\remk{\03.4}  Although an arbitrary complete variety may not have a big
line sheaf (or any nontrivial line sheaf) \cite{F, pp.~25--26 and p.~72},
a nonsingular complete variety always does.  Indeed, let $U$ be a nonempty
open affine on a nonsingular complete variety $X$, pick generators
$x_1,\dots,x_r$ for the affine ring $\Cal O_X(U)$, and let $D$ be a
Weil divisor whose support contains the polar divisors of all $x_i$.
Then $D$ is big.
\endit

\beginsection{\04}{Truncated Counting Functions}

Variations of the above conjectures using truncated counting functions can
also be made.  First of all, in Nevanlinna theory, we have:

\conj{\04.1}  Let $T$ be a Riemann surface, let $t\:T\to\Bbb C$ be a proper
surjective holomorphic map, let $X$ be a nonsingular complete complex
variety, let $K$ be the canonical
divisor class on $X$, let $\frak a\ne(0)$ be a sheaf of ideals on $X$, let $A$
be a big divisor on $X$, and let $\epsilon>0$.  Then there is a proper
Zariski-closed subset $Z$ of $X$, depending only on $\deg t$, $X$, $\frak a$,
$A$, and $\epsilon$, such that for all nonconstant holomorphic curves
$f\:T\to X$ whose images are not contained in $Z$, the inequality
$$N_f^{(1)}(\frak a,r) + N_{t,\text{Ram}}(r)
  \ge_\exc T_{K,f}(r) + T_{\frak a,f}(r) - T_{\Cal I^{-}(\frak a),f}(r)
  - \epsilon\,T_{A,f}(r) - O(1)$$
holds.
\endit

In the diophantine case, the corresponding conjecture is:

\conj{\04.2}  Let $k$ be a global field of characteristic zero, let $S$ be
a finite set of places of $k$ containing all archimedean places, let $r$ be a
positive integer, let $X$ be a nonsingular complete variety over $k$,
let $K$ be the canonical divisor class of $X$, let $\frak a\ne(0)$ be an
ideal sheaf on $X$, let $A$ be a big divisor on $X$, and let $\epsilon>0$.
Then there is a proper Zariski-closed subset $Z$ of $X$, depending only on
$k$, $S$, $X$, $\frak a$, $r$, $A$, and $\epsilon$, such that
$$N_S^{(1)}(\frak a,P) + d_k(P)
  \ge h_K(P) + h_{\frak a}(P) - h_{\Cal I^{-}(\frak a)}(P)
  - \epsilon\,h_A(P) - O(1)$$
for all $P\in (X\setminus Z)(\bar k)$ with $[k(P):k]\le r$.
\endit

It is not clear that the dependence of $Z$ on $\deg t$ or $r$ (respectively)
is necessary.

In each case, if $\frak a=\Cal O(-D)$ with $D$ an effective, reduced, normal
crossings divisor, then the above conjectures reduce to conjectures that have
already been posed; see \cite{V} for the diophantine case.

Again, we have a converse:

\prop{\04.3}  Conjectures \04.1 and \04.2 are equivalent to their respective
special cases in which $\frak a=\Cal O(-D)$ with $D$ as above.
\endit

\demo{Proof} In the diophantine case this follows by the same argument
as before.  Indeed, let $\mu\:X'\to X$, $F$, and $\eta$ be as before;
assuming that \cite{V, Conj.~2.3} holds for $F_\red$ on $X'$, we have
$$\split N_S^{(1)}(\frak a,P) + d_k(P)
  &= N_S^{(1)}(F_\red,P') + d_k(P') \\
  &\ge h_{K_{X'}+F_\red}(P') - \epsilon\,h_{\mu^{*}A}(P') - O(1) \\
  &= h_{K_X'}(P') - h_{K_{X'/X}}(P') + h_F(P')
    - h_{-K_{X'/X}+\lfloor(1-\eta)F\rfloor}(P') \\
    &\qquad- \epsilon h_A(P) - O(1) \\
  &\ge h_{K_X}(P) + h_{\frak a}(P) - h_{\Cal I^{-}(\frak a)}(P)
    - \epsilon h_A(P) - O(1)\;,
\endsplit$$
where $P'\in X'$ lies over $P\in X$.  The proof in the Nevanlinna case is
analogous.\qed
\enddemo

\comment
\beginsection{\05}{A Sample Application}

This section gives a simple example of an application of these conjectures.
Of course, since these conjectures are equivalent to earlier ones applied
on a suitable blowing-up of the variety, one can also obtain this application
by working on a blow-up.

We start with a lemma.

\lemma{\05.1}  Let $\frak a'\subseteq\frak a$ be nonzero sheaves of ideals
on $X$.  Then
$$\lambda_{\frak a}(P) - \lambda_{\Cal I^{-}(\frak a)}(P)
  \le \lambda_{\frak a'}(P) - \lambda_{\Cal I^{-}(\frak a')}(P) + O(1)
  \tag\05.1.1$$
for all $P\in X$ outside of the closed subscheme corresponding to $\frak a'$,
where $O(1)$ refers to an $M_k$\snug-constant in the diophantine case.
\endit

\demo{Proof}  Let $\mu\:X'\to X$ be a proper birational morphism satisfying
the conditions of Definition \01.1 simultaneously for $\frak a$ and $\frak a'$
(this is possible since $\mu^{*}$ is functorial).
Write $\mu^{*}\frak a=\Cal O(-F)$ and $\mu^{*}\frak a'=\Cal O(-F')$,
and pick $\epsilon>0$ such that
$\Cal I^{-}(\frak a)=\Cal I(\frak a^{1-\epsilon})$ and
$\Cal I^{-}(\frak a')=\Cal I((\frak a')^{1-\epsilon})$.

||| Not done yet.  Also mention this section in the introduction.
\enddemo
\endcomment

\enddocument